\documentclass[psamsfonts]{amsart}

\usepackage{amssymb,amsfonts}
\usepackage[all,arc]{xy}
\usepackage{enumerate}
\usepackage{mathrsfs}
\usepackage{graphicx}
\usepackage{tikz}
\usepackage{amscd}
\usepackage{amssymb}
\usepackage{xy}
\usepackage{lscape}
\usepackage{setspace}
\usepackage{amsthm}
\usepackage{mathtools}
\usetikzlibrary{matrix,arrows}

\newtheorem{thm}{Theorem}[section]
\newtheorem{cor}[thm]{Corollary}
\newtheorem{prop}[thm]{Proposition}
\newtheorem{lem}[thm]{Lemma}
\newtheorem{conj}[thm]{Conjecture}

\theoremstyle{definition}
\newtheorem{defn}[thm]{Definition}

\newtheorem{exmp}[thm]{Example}

\theoremstyle{remark}
\newtheorem{rem}[thm]{Remark}

\newcommand{\ZZ}{\mathbb Z}
\newcommand{\QQ}{\mathbb Q}

\makeatletter
\let\c@equation\c@thm
\makeatother
\numberwithin{equation}{section}

\title{Computations of the Structure of the Goldman Lie Algebra for the Torus}

\author{Felicia Tabing}

\date{}

\begin{document}

\begin{abstract}

\indent We consider the structure of the Goldman Lie algebra for the closed torus, and show that it is finitely generated over the rationals.  We also consider other traditional Lie algebra structures and determine that the Goldman Lie algebra for the torus is not nilpotent or solvable, and we compute the derived Lie algebra.\\

\end{abstract}

\maketitle

\section{The Goldman Lie Algebra}
~
\indent The Goldman Lie algebra is an algebra over the module generated by the set of free homotopy classes of loops on a surface, described by intersection and concatenation, that was introduced by William M. Goldman in 1986 \cite{Go}. \\
\indent Throughout this paper, let $\Sigma_{g,n}$ denote an oriented, genus $g$ surface with $n \geq 0$ boundary components. Denote $\hat{\pi}(\Sigma_{g,n})$ to be the set of free homotopy classes of loops on $\Sigma_{g,n}$, where we use $\hat{\pi}$ when it is clear from the context which fixed surface we are discussing. Recall the following.
\lem The set of free homotopy classes of loops on a surface $\Sigma_{g,n}$ is in one-to-one correspondence with conjugacy classes of $\pi_1(\Sigma_{g,n})$.
\rem \normalfont We can represent homotopy classes of loops by cyclically reduced words with letters the generators of the fundamental group.
\defn \normalfont Fix a surface $\Sigma_{g,n}$ and an orientation of $\Sigma_{g,n}$. Let $\alpha, \beta \in \hat{\pi}$. The \textbf{Goldman bracket} of $\alpha$ and $\beta$ is defined to be
\begin{align}
[\alpha, \beta]=\sum_{p \in \alpha \cap \beta} \epsilon (p) \alpha *_p \beta,
\end{align}
where $\alpha$ and $\beta$ intersect in transverse double points $p$, and $\epsilon (p)$ is the sign of the intersection, or $\epsilon (p)=1$ if the ordered vectors in the tangent space to $\Sigma_{g,n}$ tangent to loop $\alpha$ and $\beta$ match the orientation of the surface, and $\epsilon (p)=-1$ otherwise \cite{Go}. 
\begin{figure}
\includegraphics[scale=.4]{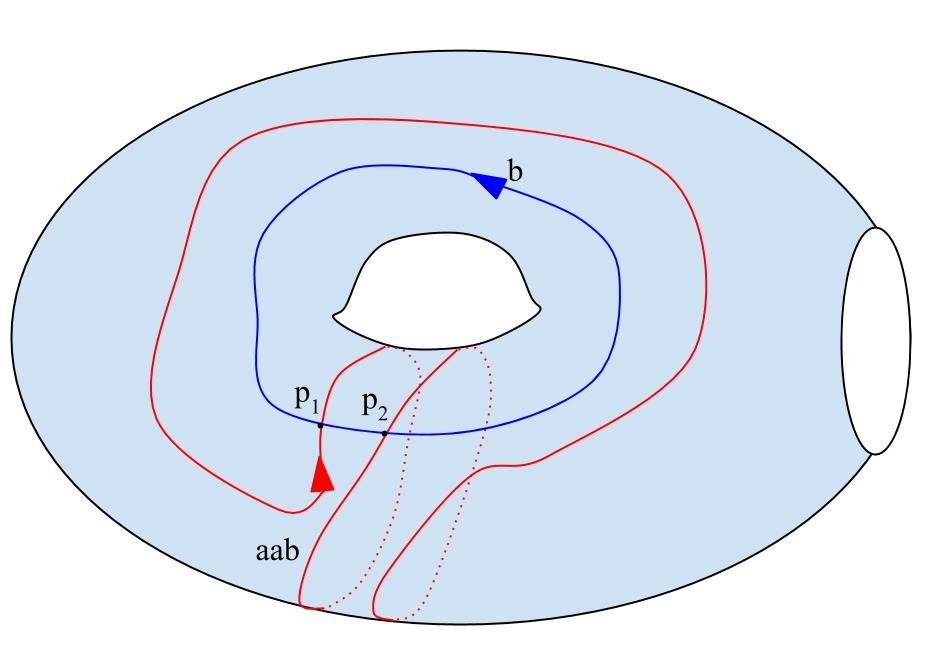}
\caption{Loops $aab$ and $b$ on the torus with one boundary component.\label{fig:sigma11}}
\end{figure}
\exmp \normalfont We compute $[aab,b]$ on the surface $\Sigma_{1,1}$. The loops represented by words $aab$ and $b$ are shown in Figure~\ref{fig:sigma11}, with intersection points $p_1$ and $p_2$. At the intersection point $p_1$, we smooth the intersection by creating a new loop, $aabb$, by following the red loop $aab$ in the direction of its orientation at $p_1$, and when returning to $p_1$, we now follow the blue loop $b$ in the direction of its orientation. When we return back to $p_1$, we close the loop. At the intersection $p_2$, we do the same, and create the loop $abab$.  We get that $[aab,b]=\pm(aabb+abab)$ where the sign depends on the chosen orientation of $\Sigma_{1,1}$.\\
\pagebreak
\thm
(Goldman) The Goldman bracket is well defined, skew-symmetric, and satisfies the Jacobi identity \cite{Go}.\\
\normalfont 
\indent We can extend the bracket linearly to $\ZZ [\hat{\pi}]$ (or $\QQ [\hat{\pi}]$), the free module over $\ZZ$ (or $\QQ$) with basis $\hat{\pi}$, to get a bilinear map
\begin{align*}
[ -,-]: \ZZ \hat{\pi}\times \ZZ \hat{\pi} \rightarrow \ZZ \hat{\pi}.
\end{align*}
Thus, $\ZZ \hat{\pi}$ is a Lie algebra with bracket $[-,-]$, which we call the \emph{Goldman Lie Algebra}, denoted by $\mathfrak{G}$ throughout the rest of this chapter. When it is unclear what the surface we are referring to, we use $\mathfrak{G}_{\Sigma_{g,n}}$
\section{Goldman Lie Algebra Structure}
~
\indent  So far, the center of the Goldman Lie algebra is known, but other parts of the structure is still unknown.  The following theorems are results on the center of the Goldman Lie algebra for closed surfaces, and surfaces with boundary, respectively.
\thm (Etingof) The center of $\mathfrak{G}_{\Sigma_{g,0}}$ is spanned by the contractible loop \cite{Et}.
\thm (Kabiraj) The center of $\mathfrak{G}_{\Sigma_{g,n}}$ is generated by peripheral loops \cite{Ka}\\
\normalfont
\indent A question posed by Chas \cite{Ch} is whether or not $\mathfrak{G}$ is finitely generated. In Goldman's paper \cite{Go}, he also introduces what is called the homological Goldman Lie algebra. This Lie algebra is defined on intersection form on the first homology group of a surface. It is known that this Lie algebra is indeed finitely generated \cite{KKT}, the ideals are known \cite{To}, and the center for the surface of infinite genus has been found \cite{KK}.\\
\section{The Goldman Lie Algebra Structure for the Torus}
\indent The closed torus is a special case as $\mathfrak{G}_{\Sigma_{1,0}}$ is finitely generated. Recall that we can represent free homotopy classes of loops on $\Sigma_{1,0}$ by cyclically reduced words in two letters, $a$ and $b$, and we can represent all homotopy classes of loops on the torus by the word $a^lb^k$ for $k,l\in \ZZ$.  The following theorem by Chas says that the Goldman bracket gives the intersection number of the free homotopy classes of loops as a coefficient of the new concatenated loop.
\prop \label{prop:torusstructure} (Chas \cite{Ch1}) The Goldman bracket structure of $\mathfrak{G}_{\Sigma_{0,1}}$ is given by
\begin{align*}
[a^ib^j,a^kb^l]=(il-jk)a^{i+k}b^{j+l}.
\end{align*}
\thm \label{thm:torus}$\mathfrak{G}_{\Sigma_{1,0}}$ is finitely generated when considered as a Lie algebra over $\QQ$.
\proof  We denote a contractible loop by $1$. We claim that $\mathfrak{G}_{\Sigma_{1,0}}$ is generated by $\{a,b,a^{-1},b^{-1}\}$. This will take many steps. We will first show that we can generate certain homotopy classes of loops. Below, we assume $n\neq 0$.
\begin{enumerate}
\item $a^nb^1=[a,a^{n-1}b]$, which we get inductively,
\item $a^n=[b^{-1},-\frac{1}{n}a^nb]$ 
\item $ab^n=[b,-ab^{n-1}]$
\item $b^n=[a^{-1},-\frac{1}{n}ab^n]$
\item $a^nb^n=[a^n,\frac{1}{n^2}b^n]$
\item $a^{-n}b=[a^{-1},-a^{-n+1}b]$
\item $a^{-n}=[b^{-1},\frac{1}{n}a^{-n}b]$
\item $a^{-1}b^n-[a^{-1},-\frac{1}{n}b^n]$
\item $a^{-n}b^n=[a^{-n},-\frac{1}{n^2}b^n]$
\item $ab^{-n}=[b^{-1},ab^{n+1}]$ which we get inductively,
\item $b^{-1}=[a^{-1},\frac{1}{n}ab^{-n}]$
\item $a^{-n}b^{-n}=[a^{-n},\frac{1}{n^2}b^{-n}]$
\item $a^nb^{-n}=[a^n,-\frac{1}{n^2}b^{-n}]$.
\item From 13. and 9. for $n=1$, we get $a^0b^0=1=[ab^{-1},\frac{1}{2}a^{-1}b]$.
\end{enumerate}
We still have a few more cases to show, namely how to generate the homotopy class of the loop $a^ib^j$ in the following cases.
\begin{enumerate}
\item[Case 1:] Suppose $i,j>0$.
\begin{enumerate}
\item Suppose $i<j$, then $j=i+r$ for some $r\in \ZZ-\{0\}$.\\ Then $\frac{1}{ar}[a^ib^i,b^j]=a^ib^j$.
\item Suppose $i>j$, then $i=j+r$ for $r\in \ZZ-\{0\}$.\\ Then $-\frac{1}{br}[a^jb^j,a^r]=a^ib^j$.
\end{enumerate}
\item[Case 2:] Suppose $i<0<j$.
\begin{enumerate}
\item Suppose $|i|<|j|$, then $j=-i+r$ for $r\in \ZZ-\{0\}$.\\ Then $\frac{1}{ar}[a^ib^{-i},b^r]=a^ib^j$.
\item Suppose $|a|>|b|$, then $i=-j+r$ for $r\in \ZZ-\{0\}$.\\ Then $-\frac{1}{br}[a^{-j}b^{j},a^r]=a^ib^j$.
\end{enumerate}
\item[Case 3:] The case $i,j<0$, and $i\neq j$ is similar to Case 1.
\item[Case 4:] The case $b<0<a$ is similar to Case 2.
\end{enumerate}
Thus, everything in $\mathfrak{G}_{\Sigma_{1,0}}$ can be generated as a Lie algebra over $\QQ$.
\qed
\cor We can refine the generators of $\mathfrak{G}_{\Sigma_{0,1}}$ to a smaller basis, namely $\{a, a^{-1}b^{-1}+b+1,b\}$.
\proof We show that we generate the basis elements mentioned in the proof of Theorem~\ref{thm:torus}.,
\begin{align*}
[a^{-1}b^{-1}+b+1,b]&=-a^{-1},\\
[a^{-1}b^{-1}+b+1,a]&=b^{-1}-ab,\\
[a,b]&=ab.
\end{align*}
\qed
\cor $\mathfrak{G}_{\Sigma_{1,0}}$ as a Lie algebra over $\QQ$ is not nilpotent, nor solvable, since $[\mathfrak{G}_{\Sigma_{1,0}} , \mathfrak{G}_{\Sigma_{1,0}}] =\mathfrak{G}_{\Sigma_{1,0}}$.
\rem $\mathfrak{G}_{\Sigma_{1,0}}$ is not finitely generated as a Lie algebra over $\mathbb{Z}$.
\proof We will show that the set $\{(n-1)a^n\}_{n>2, n\in \ZZ}$ cannot be generated. Suppose to the contrary that we can generate $(n-1)a^n$, so there exists $(i_s,j_s),(k_s,l_s) \in \ZZ ^2$ such that
\begin{align*}
\sum_{s=1}^t\pm [a^{i_s}b^{j_s},a^{k_s},b^{l_s}]=(n-1)a^n.
\end{align*}
Then, as in Proposition~\ref{prop:torusstructure},
\begin{align*}
\sum_{s=1}^t\pm [a^{i_s}b^{j_s},a^{k_s},b^{l_s}]=\sum_{s=1}^t\pm (i_sl_s-j_sk_s)a^{i_s+k_s}b^{j_s+l_s}.
\end{align*}
We need that $i_s+k_s=n$ and $j_s+l_s=0$, so $i_sl_s-j_sk_s=-i_sj_s-j_sn+j_si_s=-j_sn$. So $n\mid i_sl_s-j_sk_s$, and $n\mid \sum_{s=1}^t\pm (i_sl_s-j_sk_s)$, so $n\mid (n-1)$, which is a contradiction. \qed
\conj We conjecture that $\mathfrak{G}_{\Sigma_{g,n}}$ for $g\geq 1$ and $n>1$ is not finitely generated. For the particular case for a punctured torus, the peripheral loop is given by a commutator word. We noticed in using Chas' program for computing the bracket seems to not generate a commutator word, nor products of commutators. This needs more work, but this would mean we have a set $\{(aba^{-1}b^{-1})^n\}_{n \in \ZZ}$ of infinitely many homotopy classes of loops that each cannot be generated by any other homotopy classes of loops. 
\normalfont
\prop \label{prop:derived} The derived Lie algebra for $\mathfrak{G}_{\Sigma_{1,0}}$ is given by
\begin{align*}
[\mathfrak{G}_{\Sigma_{1,0}},\mathfrak{G}_{\Sigma_{1,0}}]=\langle d(a^ib^j),na^n,nb^n \rangle,
\end{align*}
for $d=gcd(i,j)$ and $n\in \ZZ-\{0\}$.
\proof We first show $[\mathfrak{G}_{\Sigma_{1,0}},\mathfrak{G}_{\Sigma_{1,0}}]\subset \langle d(a^ib^j),na^n,nb^n \rangle$.
\begin{enumerate}
\item[Case 1:] Suppose $d=gcd(i,j)$, $i,j \neq 0$ and $ma^ib^j \in [\mathfrak{G}_{\Sigma_{1,0}},\mathfrak{G}_{\Sigma_{1,0}}]$ for some $m\in \ZZ$. Write $xi+yj=d$ for some $x,y\in \ZZ$ and $ma^ib^j=[a^kb^l,a^p,b^q]$ for $k,l,p,q \in \ZZ$. But 
\begin{align} \label{align:gcdd}
[a^kb^l,a^p,b^q]=(kq-lp)a^{k+p}b^{l+q},
\end{align}
so we have that $k+p=i$, $l+q=j$, and $kq-lp=kj-li=d(k(\frac{j}{d})-l(\frac{i}{d}))=m$. Thus $d\mid m$.
\item[Case 2:] Suppose that $n \neq 0$ and that $[a^ib^j,a^kb^l]=ma^n$ for $i,j,k,l,m\in \ZZ$. Then $i+k=n$, $j+l=0$, so
\begin{align} \label{align:n0}
[a^ib^j,a^kb^l]=-jna^n,
\end{align} 
so $n\mid m$. 
\item[Case 3:] Showing that for $n \neq 0$ and $n\mid m$ for $mb^n \in [\mathfrak{G}_{\Sigma_{1,0}},\mathfrak{G}_{\Sigma_{1,0}}]$ is similar to Case 2.\\
\indent To show the other containment, we can consider the equality~\ref{align:gcdd} with $k=y$ and $l=-x$ for Case 1, we can consider the equality~\ref{align:n0} with $j=-1$, and we can do something similar for Case 3.
\end{enumerate}
\qed
\prop The lower central series for $\mathfrak{G}_{\Sigma_{1,0}}$ stabilizes, i.e.
\begin{align*}
[\mathfrak{G}_{\Sigma_{1,0}},G_i]=\langle d(a^ib^j),na^n,nb^n \rangle
\end{align*}
where $d=gcd(i,j)$, $n\in \ZZ-\{0\}$, for all $i\geq 0$, and $G_i=[\mathfrak{G}_{\Sigma_{1,0}},G_{i-1}]$ defined inductively, where $G_0=\mathfrak{G}_{\Sigma_{1,0}}$.
\proof For $i=1$, this is just Proposition~\ref{prop:derived}. For $i=2$, we need to show that 
\begin{align*}
[ \mathfrak{G}_{\Sigma_{1,0}}, [\mathfrak{G}_{\Sigma_{1,0}}, \mathfrak{G}_{\Sigma_{1,0}} ]]=\langle d(a^ib^j),na^n,nb^n \rangle .
\end{align*}
The "$\subset $" containment is clear. First, consider $a^ib^{-1}\in \mathfrak{G}_{\Sigma_{1,0}}$ and $a^{n-i}b \in \langle d(a^ib^j),na^n,nb^n\rangle$ (since $gcd(n-i,1)=1$). We have that 
\begin{align*}
[a^ib^{-1},a^{n-i}b]=na^n.
\end{align*}
In a similar way, we can show that $nb^n \in [ \mathfrak{G}_{\Sigma_{1,0}}, [\mathfrak{G}_{\Sigma_{1,0}}, \mathfrak{G}_{\Sigma_{1,0}} ]]$. \\
Now consider $d=gcd(i.j)$, so we can write $d=xi+yj$. Consider $a^{i+y}b^{j-x} \in \mathfrak{G}_{\Sigma_{1,0}}$ and $ a^{-y}b^x \in [\mathfrak{G}_{\Sigma_{1,0}} , \mathfrak{G}_{\Sigma_{1,0}} ]$ (since $1=\frac{i}{d}x+\frac{j}{d}y$ implies $gcd(x,y)=1$. We have
\begin{align*}
[a^{i+y}b^{j-x} , a^{-y}b^x]=da^ib^j.
\end{align*}
Thus, it follows that the lower central series stabilizes.
\qed
\cor As a Lie algebra over $\ZZ$,  $\mathfrak{G}_{\Sigma_{1,0}}$ is not nilpotent.
\normalfont

\end{document}